\documentclass[journal]{IEEEtran}
\usepackage{cite}

\ifCLASSINFOpdf
  \usepackage[pdftex]{graphicx}
\else
   \usepackage[dvips]{graphicx}
\fi
\usepackage[cmex10]{amsmath}
\usepackage{array}
\newcommand{\set}[1]{\mathbf{#1}}
\newcommand{\multskip}{}
\begin{document}

%
\title{On estimation states of hidden markov models\\ in condition of unknown transition matrix}
%
%
%


\author{Vasily~Vasilyev,
        Alexander~Dobrovidov

\thanks{V. Vasilyev is with Department of Radio Engineering and Cybernetics, Moscow Institute of Physics and Technology (State University), Moscow, Russia e-mail: evil.vasy@gmail.com}
\thanks{A. Dobrovidov is with Institute of Control Sceinces, Russian Academy of Sciences, Moscow, Russia email: dobrovidov@gmail.com}
\thanks{Manuscript received ***** **, 2015; revised ***** **, 2015.}}

%
%

\markboth{Journal of IEEE Transactions on Information Theory,~Vol.~**, No.~*, ******~2015}%
{Vasilyev \MakeLowercase{\textit{et al.}}: On estimation states of hidden markov models in condition of unknown transition matrix}
%



\maketitle

\begin{abstract}
In this paper, we develop methods of nonlinear filtering and prediction of an unobservable Markov chain with a finite set of states. This Markov chain controls coefficients of AR($p$) model. Using observations generated by AR($p$) model we have to estimate the state of Markov chain in the case of an unknown probability transition matrix. 
Comparison of proposed non-parametric algorithms with the optimal methods in the case of the known transition matrix is carried out by simulating.
\end{abstract}


%
\IEEEpeerreviewmaketitle

\section{Introduction}
%
%
%
%
\IEEEPARstart{H}{idden} markov models are very popular for modeling and simulating processes, when you do not observe...

\section{System Model}
Let $(S_n, X_n)$ be a two-component process, where $(S_n)$ is unobservable component and $(X_n)$ is observable one, $n \in \{1,2,\ldots, N\}$, $N\in\set{N}$; $(S_n)$ ``controls'' equation coefficients of $(X_n)$. Let $(S_n)$ be a stationary Markov chain with $M$ discrete states and transition matrix $\|p_{i,j}\|,\,p_{i,j}=\Pr(S_n = j \mid S_{n-1} = i)$. The process $(X_n)$ is described by the autoregressive model of order $p$:

\begin{equation}\label{observe_process}
X_n = \mu(S_n) + \sum\limits_{i=1}^p a_i(S_n)(X_{n-i} - \mu(S_n)) + b(S_n)\xi_n,
\end{equation}
where $\{\xi_n\}$ are i.i.d. random variables with the standard normal distribution, $\mu, a_i, b \in\set{R}$ are coefficients controlled by the process $(S_n)$.

As a quality measure for our methods we use mean risk $E(L(S_n, \hat{S}_n))$ with a simple loss function~$L$:

\begin{equation}\label{simple_loss_function}
L(S_n, \hat{S}_n) =	
\begin{cases}
1, &S_n \ne \hat{S}_n,\\
0, &S_n = \hat{S}_n,
\end{cases}
\end{equation}
where $\hat{S}_n = \hat{S}_n(X_1^n)$ is an estimator of $S_n$ and $X_1^n = (X_1, X_2,\ldots,X_n)$.

As known, for this risk function with the loss function~(\ref{simple_loss_function}) the optimal estimator is
\begin{equation}\label{argmax_prob}
\hat{S}_n = \underset{m\in\{1,\ldots, M\}}{\operatorname{argmax}}\Pr(S_n = m \mid X_{1}^{n}),
\end{equation}
where $\Pr(S_n = m \mid X_{1}^{n})$ is a posterior probability with respect to a $\sigma$-algebra, generated by r.v. $X_1^n$. Its realization will be denoted by
\begin{equation}
P(S_n = m \mid X_1^n = x_1^n) = P(S_n = m \mid x_1^n),
\end{equation}
where we will write $x_1^n$ instead of $X_1^n = x_1^n$.

\subsection{Basic equations}
In this paper we consider methods of filtering and prediction in the case of unknown parametres (transition matrix) of process $(S_n)$ and known parametres (equation coefficients in~\eqref{observe_process}) of process $(X_n)$. For comparison with some standard we also consider optimal filtering and prediction, where all parametres are known.

Filtering is a problem to estimate $S_n$ by using $X_1^n$. Therefore basic equations for filtering
\begin{multline}\label{filtering_basic_1}
P(S_n=m \mid x_1^n)\\
= \frac{f(x_n\mid S_n=m, x_1^{n-1})}{f(x_n \mid x_1^{n-1})}P(S_n=m\mid x_1^{n-1}),
\end{multline}
\multskip
\begin{multline}\label{filtering_basic_2}
f(x_n \mid x_1^{n-1})\\
= \sum\limits_{m = 1}^{M} f(x_n\mid S_n=m, x_1^{n-1}) P(S_n=m\mid x_1^{n-1}),
\end{multline}
can be obtained from the total probability formula. Since coefficients in~\eqref{observe_process} are known and $\xi_n \sim \mathcal{N}(0, 1)$ then 
\begin{multline}
f(x_n\mid S_n=m, x_1^{n-1})\\
= f(x_n\mid S_n=m, x_{n-p}^{n-1}) = f_{m}(x_n),\label{cond_dens_norm_1}
\end{multline}
where
\begin{multline}
f_{m}(x_n)\\
= \phi\Big(x_n; \mu(m) + \sum\limits_{i=1}^p a_i(m)(x_{n-i} - \mu(m)), b^2(m)\Big)\label{density_f_s_n}
\end{multline}
with normal probability density function
\begin{gather}\label{norm_pdf}
\phi(x;\mu,\sigma^2) = \frac{1}{\sqrt{2 \pi}\sigma}\exp\left(-\frac{(x-\mu)^2}{2\sigma^2}\right),
\end{gather}
where $x, \mu \in \set{R},.\sigma \in\set{R}^{+}$.

\section{Optimal Filtering}
In the optimal filtering all parametres are known. We use~\eqref{cond_dens_norm_1} knowing coefficients in~\eqref{observe_process} and calculate $P(S_n=m\mid x_1^{n-1})$ in~\eqref{filtering_basic_1} knowing transition matrix:
\begin{gather}
P(S_n=m\mid x_1^{n-1}) = \sum\limits_{i=1}^{M}p_{i, m} P(S_{n-1}=i \mid x_1^{n-1}). \label{markov_chain_transition}
\end{gather}
Then the~\eqref{filtering_basic_1} is transformed to the evaluation equation~\cite{dobrovidov_2012}
\begin{gather*}
P(S_n=m\mid x_1^n) = \frac{f_{m}(x_n)\sum\limits_{i=1}^{M}p_{i, m} P(S_{n-1}=i \mid x_1^{n-1})}{\sum\limits_{j = 1}^{M} f_j(x_n)\sum\limits_{i=1}^{M}p_{i, j} P(S_{n-1}=i \mid x_1^{n-1})},\label{opt_filtering_basic_1}
\end{gather*}
which will be considered as the optimal standard.

\section{Non-parametric Filtering}
\subsection{Reducing to optimization problem}
In this section, the transition matrix $\|p_{i, j}\|$ is assumed unknown, therefore we can not use the equation~\eqref{markov_chain_transition}. To overcome this uncertainty we include formula~\eqref{cond_dens_norm_1} in equations~\eqref{filtering_basic_1},~\eqref{filtering_basic_2} and obtain
\begin{gather}
P(S_n=m\mid x_1^n) = \frac{f_m(x_n)}{f(x_n \mid x_1^{n-1})}u_n(m),\label{filtering_basic_u_1}\\
f(x_n \mid x_1^{n-1}) = \sum\limits_{m = 1}^{M} f_m(x_n)u_n(m),\label{filtering_basic_u_2}
\end{gather}
where
\begin{gather*}
u_n(m) = P(S_n=m \mid x_1^{n-1}),\quad \forall m=1, \ldots, M
\end{gather*}
are new variables, which do not depend on $x_n$ and
\begin{gather*}
\sum\limits_{i=1}^M u_i = 1,\quad u_m \ge 0,\quad \forall m=1, \ldots, M.
\end{gather*}
To calculate~\eqref{filtering_basic_u_1} and~\eqref{filtering_basic_u_2} it is neccessary to find all $u_n(m)$. We need to make the assumption. We suppose that process $(S_n, X_n)$ is $\alpha$-mixing, then
\begin{gather*}
f(x_n \mid x_1^{n-1}) \approx f(x_n \mid x_{n-\tau}^{n-1}),\quad \tau\in\{1,2,\ldots,  n-1\},
\end{gather*}
and estimate density $f(x_n \mid x_{n-\tau}^{n-1})$ using kernel density estimation and designate this estimator like $\hat{f}(x_n \mid x_{n-\tau}^{n-1})$.

Let us introduce vector $\mathbf{u}_n = (u_n(1), u_n(2), \ldots, u_n(M))$ with unknown elements $u_n(m),.m=1,\ldots, M$. Then for calculating $\mathbf{u}_n$ one proposes the following estimator
\begin{multline}\label{estimator_u_n}
\hat{\mathbf{u}}_n\\
=\underset{\mathbf{u} \in \mathrm{\Delta_M}}{\operatorname{argmin\,}}\int\limits_{-\infty}^{+\infty}|\hat{f}(z_n \mid x_{n-\tau}^{n-1})- \sum\limits_{m=1}^M f_m(z_n) u_m|^2 dz_n,
\end{multline}
where
\begin{multline*}
\mathrm{\Delta}_M = \Big\{(t_1, t_2, \ldots, t_M) \in \set{R}^{M}\\
\mid \sum\limits_{i=1}^M t_i = 1, t_i \ge 0, \forall i \in\{1,2,\ldots, M\}\Big\}\
\end{multline*}
is simplex. Let us rewrite estimator $\hat{\mathbf{u}}_n$ with more detailes:
\begin{gather*}
\hat{\mathbf{u}}_n=\underset{\mathbf{u} \in \mathrm{\Delta_M}}{\operatorname{argmin\,}}I_1 - 2I_2 + I_3,
\end{gather*}
where
\begin{align*}
I_1 &= \int\limits_{-\infty}^{+\infty}\hat{f}^2(z_n \mid x_{n-\tau}^{n-1})dz_n,\\
I_2 &= \int\limits_{-\infty}^{+\infty}\sum\limits_{m=1}^M\hat{f}(z_n \mid x_{n-\tau}^{n-1})f_m(z_n) u_m dz_n,\\
I_3 &= \int\limits_{-\infty}^{+\infty}\sum\limits_{i=1}^M\sum\limits_{j=1}^M f_i(z_n) f_j(z_n) u_i u_j dz_n.
\end{align*}
Since $I_1$ does not depend on $\mathbf{u}$, then reduce it, also transform $I_2$ and $I_3$, so $\hat{\mathbf{u}}_n$ has representation
\begin{align}\label{optimization_problem_2}
\hat{\mathbf{u}}_n&=\underset{\mathbf{u} \in \mathrm{\Delta_M}}{\operatorname{argmin\,}}I_3- 2I_2\notag\\
&=\underset{\mathbf{u} \in \mathrm{\Delta_M}}{\operatorname{argmin\,}}\sum\limits_{i=1}^M\sum\limits_{j=1}^M c_{ij} u_i u_j - 2 \sum\limits_{m=1}^M c_m u_m,
\end{align}
where
\begin{align}
c_{ij} &= \int\limits_{-\infty}^{+\infty}f_i(z_n) f_j(z_n)dz_n,\label{c_ij}\\
c_{m} &= \int\limits_{-\infty}^{+\infty}\hat{f}(z_n \mid x_{n-\tau}^{n-1})f_m(z_n)dz_n. \label{c_m}
\end{align}
To solve optimization problem~\eqref{optimization_problem_2}, primarily, it is necessary to calculate latter coefficients~\eqref{c_ij} and~\eqref{c_m}, which we will obtain using kernel density estimators. Therefore we introduce following chapter.

\subsection{Kernel density estimators}
In the general case kernel density estimator of density $f$ is
\begin{gather}\label{kde_estimator}
\hat{f}(\mathbf{y; H}) = \frac{1}{N}\sum\limits_{i=1}^{N}K_{\mathbf{H}}(\mathbf{y} - \mathbf{Y}_i),
\end{gather}
where $\mathbf{y}=(y_1, y_2, \ldots, y_d)^{T}$ is argument and $\mathbf{Y}_i = (Y_{i1}, Y_{i2}, \ldots, Y_{id})^{T}$, $i=1,2, \ldots, N$ are drawn from density $f$; $K_{\mathbf{H}}(\mathbf{y}) = |\mathbf{H}|^{-1/2}K(\mathbf{H}^{-1/2}\mathbf{y})$, where $K(\mathbf{y})$ is the multivariate kernel, which is probability density function; $\mathbf{H}\in\mathcal{H}$ is the bandwidth matrix and $\mathcal{H}$ is the set of $d\times d$, symmetric and positive-definite matrixes. We propose to use unbiased cross-validation (UCV) to find $\mathbf{H}$ (univariate case proposed in~\cite{rudemo_1982},~\cite{bowman_1984} and multivariate in~\cite{sain_1994},~\cite{duong_2005}). This is a popular and relevant method is aimed to estimate
\begin{gather*}
\mathrm{ISE} (\mathbf{H}) = \int\limits_{\set{R}^d}\left(\hat{f}(\mathbf{y; H}) - f(\mathbf{y})\right)^2d\mathbf{y}
\end{gather*}
and then minimize resulting function
\begin{multline}
\mathrm{UCV}(\mathbf{H})\\
= \frac{1}{N(N-1)}\sum\limits_{i=1}^{N}\sum\limits_{\begin{smallmatrix}j = 1,\\ j\neq i\end{smallmatrix}}^{N}(K_\mathbf{H}*K_\mathbf{H} - 2K_\mathbf{H})(\mathbf{Y}_i-\mathbf{Y}_j)\\
+ \frac{1}{N}R(K)|\mathbf{H}|^{-1/2},\label{ucv}
\end{multline}
\begin{gather*}
R(K)=\int\limits_{\set{R}^d}K(\mathbf{y})^2d\mathbf{y},
\end{gather*}
where $*$ denotes a convolution. Then the estimator of $\mathbf{H}$ is
\begin{gather}\label{ucv_rule}
\mathbf{H}_{\mathrm{UCV}} = \underset{\mathbf{H} \in \mathcal{H}}{\operatorname{argmin\,}} \mathrm{UCV(\mathbf{H})}.
\end{gather}

We suppose to generate components $Y_{ik}$ of vector $\mathbf{Y}_i$ from univariate sample $x_1, x_2, \ldots, x_n$ according to the rule
\begin{gather*}
Y_{ik}= x_{(i-1)l + k},\ k=1,2,\ldots, d
\end{gather*}
where $l\in\set{N}$ influences on stochastic dependence between vectors $\mathbf{Y}_i$ (for bigger $l$ less dependence). Then we suggest to simplify obtaining of estimator~\eqref{kde_estimator} and function~\eqref{ucv}. For this aim we:

\begin{itemize}
\item use normal kernel, it means that we set equal $\mathbf{H}$ to $d$-variate normal density with zero mean vector and identity covariance matrix $\phi$;
\item use scalar $h^2$ multiple of identity $d\times d$ matrix ($\mathbf{I}_d$) for bandwidth matrix: $$\mathbf{H} = h^2 \mathbf{I}_{d}.$$
\end{itemize}

Then the estimator~\eqref{kde_estimator} becomes
\begin{multline}\label{kde_estimator_h}
\hat{f}(\mathbf{y}; h)\\
= \frac{1}{N(2\pi)^{d/2}h^d}\sum\limits_{i=1}^{N}\exp\left(-\frac{\sum\limits_{j=1}^d (y_j - x_{(i-1)l+j}) ^ 2}{2h^2}\right),
\end{multline}
with $N = 1+ \lfloor\frac{n - d}{l}\rfloor$ and the estimator of $h$ is
\begin{gather}\label{estimator_of_h}
\hat{h}=\underset{h > 0}{\operatorname{argmin\,}}\mathrm{UCV}(h),
\end{gather}
\begin{multline*}
\mathrm{UCV}(h)\\
=\frac{1}{N(N-1)(2\pi)^{d/2}h^d}\sum\limits_{i=1}^{N}\sum\limits_{\begin{smallmatrix}j = 1,\\ j\neq i\end{smallmatrix}}^{N}\frac{1}{2^{d/2}}e^{-\frac{\Delta x_{ij}}{4h^2}}-2e^{-\frac{\Delta x_{ij}}{2h^2}}\\ + \frac{1}{N(4\pi)^{d/2}h^d},
\end{multline*}
\begin{gather*}
\Delta x_{ij} = \sum\limits_{k=1}^d \left(x_{(i-1)l+k} - x_{(j-1)l+k}\right) ^ 2.
\end{gather*}
Computing minima analytically is a challenge, so a numerical calculation is popular. The function $\mathrm{UCV}(h)$ often has multiple local minima, therefore more correct way is to use brute-force search to find $\hat{h}$, however it is a very slow algorithm. In~\cite{hall_1991} it was shown that spurios local minima are more likely at too small values of $h$, so we propose to use golden section search between 0 and $h^{+}$, where
\begin{gather*}
h^{+} = \left(\frac{4}{N(d + 2)}\right)^{\frac{1}{d + 4}}\underset{k\in\{1,\ldots, d\}}{\max\,}\hat{\sigma}_k,\\
\end{gather*}
where $\hat{\sigma}_k$ is the sample standard deviation of $k$-th elements of $\mathbf{Y}_i$. The parameter $h^{+}$ is an oversmoothed bandwidth. If the matrix $\mathbf{H}$ was an unconstrained then 
\begin{gather*}
\mathbf{H}^{+}=\left(\frac{4}{N(d + 2)}\right)^{\frac{1}{d + 4}}\mathbf{S},
\end{gather*}
where $\mathbf{S}$ is a sample covariance matrix of $\mathbf{Y}_i$. The matrix $\mathbf{H}^{+}$ is oversmoothed bandwidth in the most cases. The latter estimator is proposed in~\cite{terrell_1990}. To calculate $\mathbf{H}_{\mathrm{UCV}}$ with unconstrained $\mathbf{H}$ you may use quasi-Newton minimization algorithm like in~\cite{duong_2005}.

\subsection{Calculation of coefficients $c_{ij}$ and $c_m$}
For calculating unknown coefficients $c_{ij}$ and $c_{m}$ in~\eqref{optimization_problem_2} we use formulas~\eqref{c_ij} and~\eqref{c_m}. Observe that for normal probability density function~\eqref{norm_pdf} following equation
\begin{multline*}
\int\limits_{-\infty}^{+\infty} \phi(x; \mu_1, \sigma_1^2) \phi(x; \mu_2, \sigma_2^2) dx\\
= \phi(\mu_1; \mu_2, \sigma_1 ^ 2 + \sigma_2 ^ 2) = \phi(\mu_2; \mu_1, \sigma_1 ^ 2 + \sigma_2 ^ 2)
\end{multline*}
is correct, therefore using it and~\eqref{density_f_s_n} we have
\begin{multline}\label{c_ij_finish}
c_{ij} = \int\limits_{-\infty}^{+\infty}\phi\Big(z_n; \mu(i) + \sum\limits_{k=1}^p a_k(i)(x_{n-k} - \mu(i)), b^2(i)\Big)\\
\cdot \phi\Big(z_n; \mu(j) + \sum\limits_{k=1}^p a_k(j)(x_{n-k} - \mu(j)), b^2(j)\Big)dz_n\\
=\phi\Big(\mu(i) + \sum\limits_{k=1}^p a_k(i)(x_{n-k} - \mu(i));\\
\mu(j) + \sum\limits_{k=1}^p a_k(j)(x_{n-k} - \mu(j), b^2(i) + b^2(j)\Big),
\end{multline}
also $c_{ij} = c_{j, i} > 0$. For calculating $c_{m}$ we estimate conditional density $\hat{f}(z_n \mid x_{n-\tau}^{n-1})$ applying~\eqref{kde_estimator_h}:
\begin{multline*}
\hat{f}(z_n \mid x_{n-\tau}^{n-1}) = \frac{\hat{f}(z_n, x_{n-\tau}^{n-1})}{\int\limits_{-\infty}^{+\infty} \hat{f}(z_n, x_{n-\tau}^{n-1})dz_n}\\
=\sum\limits_{i=1}^{N}\beta_{ni}(\tau)\phi(z_n; x_{(i-1)l + \tau + 1}, h^2),
\end{multline*}
\begin{gather*}
\beta_{ni}(\tau) = \frac{\exp\left(-\frac{\sum\limits_{j=-\tau}^{-1}(x_{n+j} - x_{(i-1)l+j + \tau + 1}) ^ 2}{2h^2}\right)}{\sum\limits_{k=1}^{N}\exp\left(-\frac{\sum\limits_{j=-\tau}^{-1}(x_{n+j} - x_{(k-1)l+j + \tau + 1}) ^ 2}{2h^2}\right)},
\end{gather*}
where $N = 1+ \lfloor\frac{n - 1 - d}{l}\rfloor$, bandwidth $h$ is estimated by~\eqref{estimator_of_h}. Remark that $\beta_{ni}(\tau)$ does not depend on $z_n$. Then we substitute latter estimator in~\eqref{c_m} and obtain
\begin{multline}\label{c_m_finish}
c_{m} = \int\limits_{-\infty}^{+\infty}\hat{f}(z_n \mid x_{n-\tau}^{n-1})f_m(z_n)dz_n \\
=\int\limits_{-\infty}^{+\infty}\sum\limits_{i=1}^{N}\beta_{ni}(\tau)\phi(z_n; x_{(i-1)l + \tau + 1}, h^2)\\
\cdot\phi\Big(z_n; \mu(m) + \sum\limits_{k=1}^p a_k(m)(x_{n-k} - \mu(m)), b^2(m)\Big)dz_n\\
=\sum\limits_{i=1}^{N}\beta_{ni}(\tau)\int\limits_{-\infty}^{+\infty}\phi(z_n; x_{(i-1)l + \tau + 1}, h^2)\\
\cdot\phi\Big(z_n; \mu(m) + \sum\limits_{k=1}^p a_k(m)(x_{n-k} - \mu(m)), b^2(m)\Big)dz_n\\
=\sum\limits_{i=1}^{N}\beta_{ni}(\tau)\phi\Big(x_{(i-1)l + \tau + 1};\\
\mu(m) + \sum\limits_{k=1}^p a_k(m)(x_{n-k} - \mu(m)), h^2 + b^2(m)\Big),
\end{multline}
also we remark that $c_m > 0$.

\subsection{Solution of optimization problem}
In the previous chapters we reduce main problem to optimization problem
\begin{gather*}
\hat{\mathbf{u}}_n=\underset{\mathbf{u} \in \mathrm{\Delta_M}}{\operatorname{argmin\,}}F_n(\mathbf{u}),\\
F_n(\mathbf{u})=\sum\limits_{i=1}^M\sum\limits_{j=1}^M c_{ij} u_i u_j - 2 \sum\limits_{m=1}^M c_m u_m,
\end{gather*}
where coefficients $c_{ij}$ and $c_m$ were calculated in~\eqref{c_ij_finish} and~\eqref{c_m_finish}. Let us consider kind of optimization. We have that $\mathrm{\Delta}_m$ is convex set and Hessian matrix of function $F_n(\mathbf{s})$ is
\begin{gather*}
\mathcal{L}''_{\mathbf{u}} = 2\cdot\begin{pmatrix}
c_{11} & c_{12} & \ldots & c_{1M}\\
c_{21} & c_{22} & \ldots & c_{2M}\\
\vdots & \vdots & \ddots & \vdots\\
c_{M1} & c_{M2} & \ldots & c_{MM}
\end{pmatrix}.
\end{gather*}
If $\mathcal{L}''_{\mathbf{u}}$ is positive defined, then $F_n(\mathbf{s})$ is convex, thus we have convex optimization. In this case we propose to use Karush--Kuhn--Tucker (KKT) conditions~\cite{kuhn_1951},~\cite{boyd_2004}, because of:
\begin{itemize}
\item our case is special because there is opportunity to solve KKT conditions analytically;
\item for convex optimization KKT conditions, which are primarily necessary, are also sufficient;
\end{itemize}
else you may apply methods of quadratic programming. Also we want remark that $\mathcal{L}''_{\mathbf{u}}$ does not depend on variables $u_i$ and coefficients $c_m$, which means that previous kernel density estiamtors have no influence on kind of optimization.

Let us consider KKT conditions, then Lagrangian is
\begin{gather*}
\mathcal{L} = \lambda_0 F_n(\mathbf{u}) + \sum\limits_{i=1}^{M}\lambda_{i} (-u_i) +  \lambda_{M+1}\left(\sum\limits_{i=1}^M u_i - 1\right),
\end{gather*}
where $\lambda^* = (\lambda_0^*, \lambda_1^*, \ldots, \lambda_{M+1}^*)\in\set{R}^{M+2}$. We need to find $\lambda^*$ and $\mathbf{u}^*$ such that stationary condition
\begin{gather*}
\mathcal{L}'_{u_i} = 2\lambda_0^*\left(\sum\limits_{j=1}^M c_{ij} u_j^* - c_i\right) - \lambda_{i}^* + \lambda_{M+1}^* = 0,\\
\forall i=1, \ldots, M
\end{gather*}
primal feasibility
\begin{gather*}
-u_i^* \le 0,\quad \forall i=1, \ldots, M\,\\
\sum\limits_{i=1}^M u_i^* - 1 = 0,
\end{gather*}
dual feasibility
\begin{gather*}
\lambda_{i}^* \ge 0,\quad \forall i=1, \ldots, M\
\end{gather*}
complementary slackness
\begin{gather*}
\lambda_{i}^* u_i^* = 0,\quad \forall i=1, \ldots, M
\end{gather*}
hold. Let $\lambda_0^*=0$ to check that the gradients of constraints are linearly independent at $\mathbf{u}^*$, so KKT conditions lead to system
\begin{gather*}
\begin{cases}
\lambda_1^* = \lambda_{2}^*=\ldots=\lambda_{M+1}^*,\\
\lambda_{i}^* u_i^* = 0,\ \lambda_{i}^* \ge 0, &\forall i=1, \ldots, M\\
\sum\limits_{i=1}^M u_i^* = 1,\ u_i^* \ge 0, &\forall i=1, \ldots, M
\end{cases}
\end{gather*}
which could be solved only with $\lambda^* = \vec{0}$, which means that gradients of constraints are linearly independent for any $\mathbf{u}^*$.
The vector $\lambda^*$ is defined with an accuracy of $\alpha > 0$, so we define $\lambda_0=1/2$, then KKT conditions lead to a system
\begin{gather*}
\mathbf{C}\cdot\vec{\rho} = \mathbf{c},
\end{gather*}
where
\begin{gather*}
\mathbf{C} = \begin{pmatrix}
c_{11} & c_{12} & \cdots & c_{1M} & -1 & 0 &\cdots &0 & 1\\
c_{21} & c_{22} & \cdots & c_{2M} & 0 & -1 & \cdots &0 & 1\\
\vdots & \vdots & \ddots & \vdots & \vdots& \vdots & \ddots & \vdots & \vdots\\
c_{M1} & c_{M2} & \cdots & c_{MM} & 0 & 0&\cdots & -1 & 1\\
1 & 1 & \cdots & 1 & 0 & 0 & \cdots & 0 & 0
\end{pmatrix},\\
\vec{\rho} = \begin{pmatrix}
u_1^*\\
\vdots\\
u_M^*\\
\lambda_1^*\\
\vdots\\
\lambda_{M+1}^*
\end{pmatrix},\ \mathbf{c} = \begin{pmatrix}
c_1\\
c_2\\
\vdots\\
c_M\\
1
\end{pmatrix},\\
\lambda_{i}^* u_i^* = 0,\ \lambda_{i}^* \ge 0,\ u_{i}^* \ge 0.\quad \forall i=1, \ldots, M
\end{gather*}
To solve last system it is necessary to consider all combinations of pairs $(u_i^*, \lambda_i^*),\ \forall i=1, \ldots, M$, where $u_i^*$ or $\lambda_i^*$ is equal to 0 (not both). Total amount of combinations is equal to $2^M$. If $u_i^* = 0$ then $i$-th column in the matrix $\mathbf{C}$ and $i$-th row in $\bar{\rho}$ are reduced, else $\lambda_i^*=0$ and $(M + i)$-th column in the matrix $\mathbf{C}$ and $(M + i)$-th row in $\bar{\rho}$ are reduced. After choosing zero element in each pair $(u_i^*, \lambda_i^*),\ \forall i=1, \ldots, M$ matrix $\mathbf{C}$ is reduced to an $(M + 1)\times (M + 1)$-matrix $\mathbf{C}_r$ and $\bar{\rho}$ to $(M + 1)\times 1$-matrix $\bar{\rho}_r$. Therefore for each combination it is necessary to calculate
\begin{gather*}
\bar{\rho}_r = \mathbf{C}_r^{-1}\cdot\mathbf{c}.
\end{gather*}
If the first $M$ elements in $\bar{\rho}_r$ are non-negative then obtained $\mathbf{u}^*$ is a solution ($\hat{\mathbf{u}}_n$) of optimization problem and there is no reason to calculate $\bar{\rho}_r$ for the next combination, because in convex optimization local minima is global minima.

As a result, we substitute estimator $\hat{\mathbf{u}}_n$ in~\eqref{filtering_basic_u_1} and~\eqref{filtering_basic_u_2} and problem of non-parametric filtering is solved.

\section{One-step Ahead Prediction}
We will consider one-step ahead prediction. Like for filtering we minimize mean risk $E(L(S_n, \hat{S}_n))$ with simple loss function~\eqref{simple_loss_function}. Therefore optimal estimator of $S_n$ is
\begin{gather*}
\hat{S}_n = \underset{m\in\{1,\ldots, M\}}{\operatorname{argmax}}\Pr(S_n = m \mid X_{1}^{n-1}).
\end{gather*}
We remark that probabilty $\Pr(S_n = m \mid X_{1}^{n-1})$ is already obtained in the considered approaches of filtering: for optimal prediction it is written in~\eqref{markov_chain_transition} and for non-parametric prediction accordingly in~\eqref{optimization_problem_2}. It means that we primarily solve problem of one-step ahead prediction and then filtering problem.

\section{Example}
\begin{figure}[!t]
\centering
\includegraphics[scale=1]{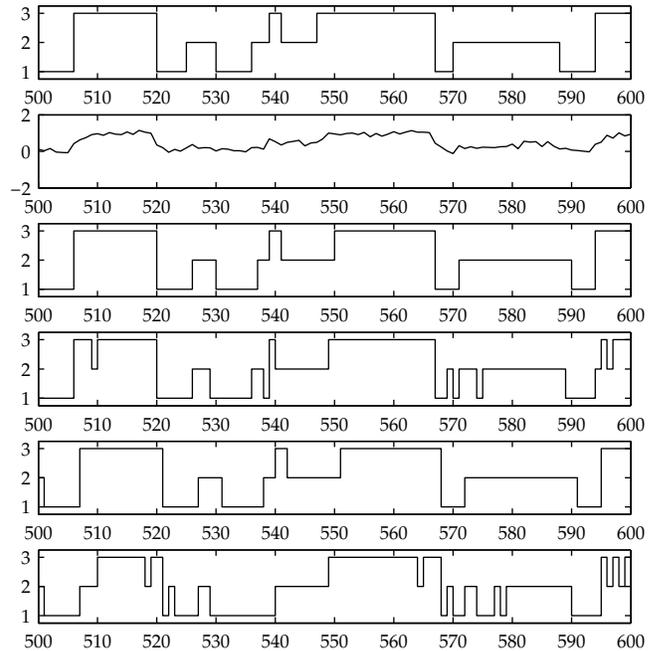}
\caption{From top to bottom: 1 --- unobservable $s_n$; 2 --- observable $x_n$;\quad\quad\quad 3, 4 --- optimal and non-parametric filtering; 5, 6 --- optimal and non-parametric prediction.}
\label{example_figures}
\end{figure}

Let the Markov chain $(S_n)$ has 3 states ($M=3$) and transition matrix
\begin{gather}
\|p_{i, j}\| = \begin{pmatrix}
0.8 & 0.1 & 0.1\\
0.05 & 0.9 & 0.05\\
0.1 & 0.05 & 0.85
\end{pmatrix}.
\end{gather}
Sample volume $n$ is changed from 500 to 600. Observable process $(X_n)$ is simulated like AR($2$) model with coefficients $\mu\in\{0, 0.5, 1\}$, $a_1\in\{0.3, 0.2, 0.1\}$, $a_2\in\{0.2, 0.3, 0.4\}$, $b\in\{0.1, 0.2, 0.1\}$. Also we take $\tau=2$ and $l=1$. The results are presented in \figurename~\ref{example_figures} and sample mean errors after 50 repeated experiments in Table~\ref{example_errors}.

\begin{table}[!t]
\renewcommand{\arraystretch}{1.3}
\caption{Sample Mean Errors}
\label{example_errors}
\centering
\begin{tabular}{ccc}
\hline
& Filtering error, \% & Prediction error, \%\\
\hline
Optimal & 16.4 & 26.6\\ 
Non-parametric & 22.7 & 37.6\\
\hline
\end{tabular}
\end{table}

\section{Conclusion}
Preparing...


%





\ifCLASSOPTIONcaptionsoff
  \newpage
\fi



\bibliographystyle{IEEEtran}
\bibliography{IEEEabrv,mybib}

\begin{thebibliography}{1}
\providecommand{\url}[1]{#1}
\csname url@samestyle\endcsname
\providecommand{\newblock}{\relax}
\providecommand{\bibinfo}[2]{#2}
\providecommand{\BIBentrySTDinterwordspacing}{\spaceskip=0pt\relax}
\providecommand{\BIBentryALTinterwordstretchfactor}{4}
\providecommand{\BIBentryALTinterwordspacing}{\spaceskip=\fontdimen2\font plus
\BIBentryALTinterwordstretchfactor\fontdimen3\font minus
  \fontdimen4\font\relax}
\providecommand{\BIBforeignlanguage}[2]{{%
\expandafter\ifx\csname l@#1\endcsname\relax
\typeout{** WARNING: IEEEtran.bst: No hyphenation pattern has been}%
\typeout{** loaded for the language `#1'. Using the pattern for}%
\typeout{** the default language instead.}%
\else
\language=\csname l@#1\endcsname
\fi
#2}}
\providecommand{\BIBdecl}{\relax}
\BIBdecl

\bibitem{dobrovidov_2012}
A.~V. Dobrovidov, G.~M. Koshkin, and V.~A. Vasiliev, \emph{Non-parametric
  models and statistical inference from dependent observations}.\hskip 1em plus
  0.5em minus 0.4em\relax USA: Kendrick Press, 2012.

\bibitem{rudemo_1982}
M.~Rudemo, ``Empirical choice of histograms and kernel density estimators,''
  \emph{Scandinavian Journal of Statistics}, vol.~9, pp. 65--78, 1982.

\bibitem{bowman_1984}
A.~Bowman, ``An alternative method of cross-validation for the smoothing of
  density estimates,'' \emph{Biometrika}, vol.~7, pp. 353--360, 1984.

\bibitem{sain_1994}
S.~R. Sain, K.~A. Baggerly, and D.~W. Scott, ``Cross-validation of multivariate
  densities,'' \emph{Journal of the American Statistical Association}, vol.~89,
  pp. 807--817, 1994.

\bibitem{duong_2005}
T.~Duong and M.~L. Hazelton, ``Cross-validation bandwidth matrices for
  multivariate kernel density estimation,'' \emph{Scandinavian Journal of
  Statistics}, vol.~32, no.~3, pp. 485--506, 2005.

\bibitem{hall_1991}
P.~Hall and J.~Marron, ``Local minima in cross-validation functions,''
  \emph{Journal of the Royal Statistical Society, Series B (Methodological)},
  vol.~53, pp. 245--252, 1991.

\bibitem{terrell_1990}
G.~R. Terrell, ``The maximal smoothing principle in density estimation,''
  \emph{J. Amer. Statist. Assoc.}, vol.~85, pp. 470--477, 1990.

\bibitem{kuhn_1951}
H.~W. Kuhn and A.~W. Tucker, ``Nonlinear programming,'' in \emph{Proceedings of
  the Second Berkeley Symposium on Mathematical Statistics and Probability},
  Berkeley, California, 1951, pp. 481--492.

\bibitem{boyd_2004}
S.~Boyd and L.~Vandenberghe, \emph{Convex Optimization}.\hskip 1em plus 0.5em
  minus 0.4em\relax Cambridge: Cambridge University Press, 2004.

\end{thebibliography}
%



%

\begin{IEEEbiography}{Vasily Vasilyev}
Preparing...
\end{IEEEbiography}




\begin{IEEEbiography}{Alexander Dobrovidov}
Preparing...
\end{IEEEbiography}




\end{document}